\documentclass[12pt,a4paper]{amsart}
\usepackage{amsmath}
\usepackage{amssymb, amsthm, amscd}
\usepackage[english]{babel}

\advance\hoffset-20mm \advance\textwidth40mm
\newtheorem{theorem}{Theorem}
\newtheorem{lemma}{Lemma}

\newtheorem{definition}{Definition}
\newtheorem{remark}{Remark}

\def\Real{{\mathbb R}}
\def\Compl{{\mathbb C}}
\def\Affspace{{\mathbb A}}

\def\Aut{\mathop{\rm Aut}}
\def\SAut{\mathop{\rm SAut}}
\def\Id{\mathop{\rm Id}}

\def\Stab{\mathop{\rm Stab}}
\def\Spec{\mathop{\rm Spec}}
\def\Der{\mathop{\rm Der}}

\def\exp{\mathop{\rm exp}}
\def\reg{\mathop{\rm reg}}
\def\rk{\mathop{\rm rk}}

\def\Int{\mathop{\rm Int}}
\def\SP{\mathop{\rm Susp}}

\def\ll1{l_{\lambda}^{-1}(1)}
\def\lm1{l_{\mu}^{-1}(1)}

\begin{document}

\title{Infinitely transitive actions on real affine suspensions}

\author[au1]{Karine Kuyumzhiyan}
\address[au1]{
National Research University Higher School of Economics, Laboratory of Algebraic Geometry and its Applications, Russia; IUM, Moscow, Russia; \\ 
Laboratoire J.-V.Poncelet (UMI 2615 of CNRS), Moscow, Russia} 
\email{karina@mccme.ru}
\thanks{The first author was supported by AG Laboratory HSE, RF government grant, ag.
11.G34.31.0023, by the "EADS Foundation Chair in Mathematics", Russian-French Poncelet Laboratory (UMI 2615 of CNRS), and Dmitry Zimin fund "Dynasty".}
\author[au2]{Fr\'ed\'eric Mangolte}
\address[au2]{Universit\'e d'Angers (LAREMA),
UMR 6093 et FR 2962 du CNRS, France}
\email{Frederic.Mangolte@univ-angers.fr}

\begin{abstract} 
A group $G$ acts infinitely transitively on a set $Y$ if for every positive integer~$m$, its action is $m$-transitive on~$Y$. Given a real affine algebraic variety $Y$ of dimension greater than or equal to~$2$, we show that, under a mild restriction, if the special automorphism group of $Y$ (the group generated by one-parameter unipotent subgroups) is infinitely transitive on each connected component of the smooth locus $Y_{\reg}$, then for any real affine suspension~$X$ over~$Y$, the special automorphism group of~$X$ is infinitely transitive on each connected component of $X_{\reg}$. This generalizes a recent result by Arzhantsev, Kuyumzhiyan, and Zaidenberg over the field of real numbers. 
\end{abstract}
%

\maketitle

\section*{Introduction}
In this note the algebraic varieties are affine with ground field of characteristic zero. Let~$Y$ be such a variety and let $f$ be a non-constant polynomial function on~$Y$. Recall that the \textit{suspension} over $Y$ along $f$ is the hypersurface $X\subseteq Y\times \Affspace^2$ given by the equation
$uv-f(y)=0$.

The paper~\cite{AKZ} shows that in many situations the infinite homogeneity of an affine variety  induces the infinite homogeneity of its iterated suspensions.
Namely, if the special automorphism group of an affine variety $Y$ of dimension at least two  acts infinitely transitively on the smooth locus of $Y$ and if at every smooth point of $Y$ 
the tangent space is spanned by the tangent vectors to the orbits of
one-parameter additive subgroups, then every suspension over~$Y$ satisfies the same two properties. The proof in such generality is however valid provided that the ground field is algebraically closed. When the ground field is~$\Real$, it is proved that the same result holds under two restrictions: the smooth locus of $Y$ is connected and the function $f$ is surjective.
The aim of this note is to settle
 the real case for any $Y$ and any $f$.

Note that the notion of affine suspension was introduced in \cite{KZ} as a particular instance of an affine modification. The latter is an important tool to understand the structure of birational morphisms between affine varieties, see~\cite{D05}. Note also that infinite transitivity, sometimes called very transitivity, was recently studied in the context of real algebraic geometry, see~\cite{hm3,hm4,bm1}.

\section{Infinitely transitive actions}

We recall the notations and state the main results. 

\begin{definition}\label{def3}
A  \textit{suspension} over an affine variety $Y$ is a hypersurface $X\subseteq Y\times \Affspace^2$ given by an equation
$uv-f(y)=0$, where $f\in\Real[Y]$ is non-constant. 
In particular, $\dim X=1+\dim Y$, and the projection on the first factor induces a natural map 
$$
\pi\colon X\to Y.
$$
\end{definition}

\begin{definition}
We say that the action of a group $G$ on a set $Y=Y^1\sqcup Y^2\sqcup \ldots \sqcup Y^s$ is {\it infinitely transitive on each connected component} if for every
$s$-tuple $(m_1,\ldots,m_s)$, it is transitive on $(m_1+\ldots+m_s)$-tuples of the form 
$$
(P_1^{1},\ldots, P^1_{m_1}, P^2_1,\ldots,P^2_{m_2},\ldots,P^s_1,\ldots,P^s_{m_s})\,,
$$ 
where $P^i_j\in Y^i$ are pairwise distinct.
\end{definition}

For an algebraic variety $X$, let the special automorphism group $\SAut(X)$ be the subgroup of $\Aut(X)$ generated by all of its one-parameter
subgroups isomorphic to the additive group~$(\Real,+)$.
Note that the action of $\SAut(X)$ does not mix regular and singular points. 

Let $Y$ be an algebraic variety over $\Real$. 
We say that a point
$y\in Y$ is \textit{flexible} if the tangent space $T_y Y$ is
spanned by the tangent vectors to the orbits $H \cdot y$ of
one-parameter subgroups~$H\subseteq\SAut (Y)$, $H\cong (\Real,+)$. The
variety $Y$ is called \textit{flexible} if every smooth point~$y\in Y_{\reg}$ is.    

\begin{theorem}[Infinite transitivity on each connected component]\label{main}
Let $Y$ be an affine algebraic variety defined over $\Real$ and $f\in \Real[Y]$. Assume that for each connected component $Y^i$ of $Y_{\reg}$, the dimension $\dim Y^i \geq 2$ and $f$ is non-constant on~$Y^i$. 

If $Y$ is flexible and the action of $\SAut(Y)$ on $Y_{\reg}$ is infinitely transitive on each connected component, then the suspension $X=\SP(Y,f)$ is flexible and 
$\SAut(X)$ acts on $X_{\reg}$ infinitely transitively on each connected component.

\end{theorem}

\begin{remark}
When $Y_{\reg}$ is connected and $f(Y_{\reg})=\Real$, the result is given by~\cite[Theorem~3.3]{AKZ}.
Notice that under these conditions, $X_{\reg}$ is also connected.
\end{remark}

\begin{remark}If $X_{\reg}$ is not connected, then $\SAut(X)$ is not even one-transitive on~$X_{\reg}$. 
Indeed, the action of $\SAut(X)$ on $X$ fixes each connected component of~$X$: every special automorphism  $g$ admits a decomposition $\prod h_j(1)$, where each~$h_i$ is a one-parameter additive group. For any $x\in X$, the arc  $t \to \prod h_i(t) \cdot x$ then connects $x$ to $g\cdot x$. \end{remark}

\begin{remark}
The number of connected components can grow on each step even if we started with a variety $Y$ whose non-singular part $Y_{\reg}$ is connected. 
Indeed, if $f$ is positive on one of the connected components, say on $Y^1$, then the set $\{(y,u,v) \,|\, uv=f(y), y\in Y^1\}$ 
splits into $\{(y,u,v),u>0,v>0\}$ and $\{(y,u,v),u<0,v<0\}$. We may choose one connected component of the suspension and further perform suspensions over this connected component. 
\end{remark}

As a preliminary part of Theorem~\ref{main}, we prove the following theorem.

\begin{theorem}[Infinite transitivity on one connected component]\label{vtcc}
Let $Y$ be an affine algebraic variety defined over $\Real$ and $f\in \Real[Y]$. Assume that $Y_{\reg}$ contains a flexible connected component $Y^1$ of dimension at least $2$ such that $\SAut(Y)$ acts infinitely transitively on $Y^1$ 	and~$f|_{Y^1}$ is non-constant. Let $X^1$ be a connected component of the smooth locus of the suspension~$\SP(Y^1,f)\subseteq X=\SP(Y,f)$. Then $X^1$ is flexible and 
$\SAut(X)$ acts infinitely transitively on $X^1$.
\end{theorem}

Note that the new paper \cite{AFKKZ} proves that over an algebraically closed field, an affine variety $X$ of dimension $\geq 2$ is flexible if and only if  $\SAut(X)$ is infinitely transitive on $X_{\reg}$ and, even more, if and only if  $\SAut(X)$ is one-transitive.

\section{Affine modifications and lifts of automorphisms}

In this section we prove the basic results of the theory over the real numbers. The main part is close to the treatment in \cite[\S3]{AKZ}.
\medskip

For every  geometrically irreducible algebraic variety $X$ over the ground field~$\Compl$, there is a natural one-to-one correspondence between locally nilpotent derivations (LND's)  $\delta$ 
on~$\Compl[X]$ and algebraic actions of one-parameter subgroups $ (\Compl,+)\cong H_{\delta}\subseteq \SAut(X)$. Namely, given a locally nilpotent derivation $\delta$, the corresponding action is the exponential 
$
(t,f) \mapsto \sum_{k=0}^\infty \frac{t^k}{k!}\delta^k(f)
$.
Conversely, for every algebraic action $\sigma$ of a subgroup $(\Compl,+)\cong H\subseteq \SAut(X)$, the derivation along the tangent vector field to the orbits of $\sigma$, given by $\frac{\sigma(t,f)-f}{t}\vert_{t=0}$ is an LND, see \cite[\S~1.5]{freu}.
The lemma below shows that the same is true for $\Real$. 

Let $G$ be a group. Recall that a $G$-module $V$ is {\it rational} if each $v\in V$ belongs to a finite dimensional $G$-invariant linear subspace $W\subseteq V$ and the $G$-action on~$W$ defines a homomorphism of algebraic groups $G\to GL(W)$. A $G$-algebra is an algebra with a structure of $G$-module.

\begin{lemma}
There are one-to-one correspondences between locally nilpotent derivations of~$\Real[X]$, unipotent subgroups $(\Real,+)\subseteq \Aut(X)$, and structures of rational $(\Real,+)$-algebras on $\Real[X]$.
\end{lemma}
\begin{proof}

For an LND $D$, the corresponding  $(\Real,+)$-algebra is defined  by the following formula:
$$
t\colon f = \exp (tD)(f)=f+tD(f)+\frac{t^2}{2!}D^2(f)+\ldots.
$$

For any fixed $f$ there exists a natural $N$ with $D^N(f)=0$, so this formula gives a polynomial in~$t$. Hence, $f$ belongs to a $(\Real,+)$-invariant linear subspace $\langle f, D(f),\ldots,D^{N-1}(f) \rangle$, which shows that $\Real[X]$ is rational as a $(\Real,+)$-algebra.

Conversely, let  $(A, t \mapsto \varphi_t)$ be a rational $(\Real,+)$-algebra. Let us define

$$
D(f)=\frac{d}{dt}\mid_{t=0}\varphi_t(f),\; f\in A.
$$

The main point is to prove that for each $f\in A$ some power $D^N(f)$ vanishes.
Consider a finite dimensional invariant subspace $W\subseteq \Real[X]$, $f\in W$. Obviously,~$D$ preserves $W$ and the action of $\exp(D)$ on $W\otimes_\Real \Compl$ is unipotent. By the Lie-Kolchin theorem, the action of~$D$ is upper-triangular  in some basis of 
$W\otimes _\Real \Compl$. This means in particular that $D^N=0$ for some~$N$. Note that the actions of $D$ on $W$ and of $D$ on $W\otimes _\Real \Compl$ were originally given by the same matrix, hence, this matrix is nilpotent,
and the derivation $D$ on $\Real[X]$ is in fact locally nilpotent. 
\end{proof}

Here is the geometric counterpart of the affine suspensions introduced above. Let $X=\SP(Y,f)$ be a suspension of $Y$ given by Definition~$\ref{def3}$. Consider the cylinder $Y\times \Affspace^1$ over $Y$,
where $\Affspace^1=\Real[v]$. Then $\SP(Y,f)$ is the blow-up of $Y\times \Affspace^1$ with center $(f,v)$ along $v$, which is a particular instance of an affine modification, see~\cite[Example 1.4]{KZ}.

Let $\delta_0$ be an LND on $\Real[Y]$ and let $H_{\delta_0}$ be the associated $(\Real,+)$-action on $Y$. 
Recall the construction of an LND $\delta_1$ which lifts $\delta_0$ to  $X$ (see \cite[Lemma~3.3]{AKZ} or \cite{KZ}). Let $\delta'$ be the lift of $\delta_0$ 
on $\Real[Y\times \Affspace^1]$ defined by $\delta'(v)=0$ and consider a product $\delta_1 = q\delta'$ by a polynomial $q(v)$ such that $q(0)=0$. 
Choosing $q$ such that the value of $\delta_1$ on $u$ preserves the relation $\delta_1(uv-f(y))=0$, we get an LND on $\Real[X]$ which satisfies
 
\begin{eqnarray}\label{drob}
\delta_1(g)&=&q(v)\delta_0(g)\quad \mbox{ for all functions } g\in\Real[Y],\nonumber\\
\delta_1(v)&=&0,\\
\delta_1(u)&=& \frac{q(v)}{v}\,\delta_0(f).\nonumber
\end{eqnarray}
There is some freedom in the choice of $q(v)$. All the derivations obtained in this way annihilate the function $v\in\Real[X]$, and the corresponding actions preserve the sections $V_c=\{v=c\}\cap X$. Notice that $X$ can also be considered as the blow-up of $Y\times \Spec\Real[u]$, and the lifts of LNDs obtained in this way annihilate the function $u\in \Real[X]$.

We denote by $G_v$ (resp. $G_u$) the subgroup of $\SAut(X)$ generated by one-parameter subgroup lifted from $Y\times \Spec\Real[v]$ (resp. $Y\times \Spec\Real[u]$).
Recall the following.
\begin{lemma}
\cite[Lemma~3.2]{AKZ}
\label{regular} 
Let $Y$ be an affine variety over a field of characteristic $0$ and $X$ be a suspension of $Y$. Then the restriction $\pi \colon X\subset Y\times\Affspace^2\to Y$ of the canonical projection satisfies
$
\pi (X_{\reg}) = Y_{\reg}
$.
\end{lemma}

We denote by $Y_{\reg}=Y^1\sqcup\ldots\sqcup Y^s$ the decomposition of $Y_{\reg}$ into connected components. If $f$ is not surjective, then the suspension over a connected component~$Y^i$ of $Y_{\reg}$ is either connected or consists of two components: if $f|_{Y^i}$ does not attain zero, it may be assumed positive, $u$ and $v$ neither attain zero, but can be either both negative, or both positive. 

For every $c\in\Real$ the hyperplane section $\{v=c\}\subset X$ will be denoted by $V_c$. We denote by~$v(P)$ the $v$-coordinate of a point $P\in X$.

Given $k$ distinct constants $c_1,\ldots,c_k\in\Real\setminus \{0\}$, we let $\Stab^v_{c_1\ldots c_k}$ be the subgroup of
$G_v$ fixing pointwise the hypersurfaces $V_{c_s}\subseteq X$, $s=1,\ldots, k$.
Observe that, as a subgroup of~$G_v$, the group $\Stab^v_{c_1\ldots c_k}$ stabilizes all the levels $V_c$ of the function $v\in \Real[X]^{G_v}$. Likewise, let $\Stab^u_{c_1\ldots c_k}\subset G_u$ be the subgroup of maps inducing the identity on the levels $U_{c_s}$ of the function $u\in \Real[X]^{G_u}$.

\begin{lemma}\label{stablevel} 
If the action of $\SAut(Y)$ on $Y_{\reg}$ is infinitely transitive on each connected component,
then for every  distinct values  $c_0,c_1,\ldots,c_k\in\Real\setminus \{0\}$, the
group $\Stab^v_{c_1\ldots c_k}$  acts infinitely transitively on each connected component of
$V_{c_0}\cap X_{\reg}$. The same is true for the action of $\Stab^u_{c_1\ldots c_k}$ 
on~$U_{c_0}\cap X_{\reg}$.
\end{lemma}

\begin{proof}
Let $P_1,\ldots,P_m$ and $Q_1,\ldots,Q_m$ be two $m$-tuples of distinct points of $V_{c_0}\cap X_{\reg}$. Since~$\pi$ restricts to an isomorphism $V_{c_0}\cap X_{\reg}\to \pi(X_{\reg})=Y_{\reg}$, we get $\pi(P_j)\ne\pi(P_l)$ and $\pi(Q_j)\ne\pi(Q_l)$ for $j\ne l$. Moreover, two points belong to the same connected component of $V_{c_0}\cap X_{\reg}$ if and only if their projections belong to the same connected component of~$Y_{\reg}$. As a consequence, there exists a special automorphism $\psi$ such that $\psi\cdot \pi(P_j)=\pi(Q_j), \forall j$. The special automorphism $\psi$ decomposes into exponentials of LND's. 
We lift each of these derivations using the polynomial $q(z)=\alpha z(z-c_1)\dots(z-c_k)$ where $\alpha\in\Real\setminus \{0\}$ is determined by $q(c_0)=1$. (Compare~\cite[Lemma 3.4]{AKZ}.)
\end{proof}

\begin{lemma}\label{level} 
Let $Y^1\subset Y$ be a connected component of $Y_{\reg}$ of dimension at least two. Then for every continuous
function $f\colon Y \to\Real$ and for each $c\in\Int f(Y^1)$, the level set $f^{-1}(c)$ is infinite.
\end{lemma}

\begin{proof}(See~\cite[Lemma 3.6]{AKZ})
We may assume that $f|_{Y^1}$ is non-constant. Choose two points $y_1,y_2\in Y^1$ such that $f(y_1)=c_1<c$ and ${f(y_2)=c_2>c}$. They can be
joined by a smooth path $l$ in $Y^1$. There exists a tubular neighborhood $U$ of $l$
diffeomorphic to a cylinder~$\Delta\times I$, where $I=[0,1]$ and $\Delta$ is a ball of dimension $\dim \Delta=\dim Y^1-1\ge 1$. So
there is a continuous family of paths joining $y_1$ and $y_2$ within $U$ such that any two of them meet only at their ends $y_1$ and
$y_2$. By continuity, on each of these paths there is a point $y'$ with $f(y')=c$. In particular, the level set $f^{-1}(c)$ is infinite.
\end{proof}

\begin{lemma}\label{avoid0} Let $Y$ be an affine variety over $\Real$ and $X$ be a suspension of $Y$.
Let $Y^1\subset Y$ be a connected component of $Y_{\reg}$ of dimension at least two, $f\in \Real[Y]$ such that $0\in \Int f(Y^1)$, and $X^1$ be the suspension over $Y^1$. 
If $Y^1$ is flexible and the action of $\SAut(Y)$ is infinitely transitive on~$Y^1$, then for every set of distinct points $P_1,\ldots,P_m$ of $X^1$ there exists a special automorphism 
$\varphi \in \SAut(X)$ such that $\varphi\cdot P_j\not\in U_0\cup V_0$ for all~$j$.
\end{lemma}

\begin{proof} 
We follow the proof of~\cite[Lemma 3.5]{AKZ}.
We say that the point~$P_i=(R_i, u_i,v_i)\in X^1$ is {\it hyperbolic} if $u_iv_i\ne 0$, i.e. $P_i\not\in U_0\cup V_0$. We
have to show that the original collection can be moved by means of a special automorphism so that all the points become hyperbolic.
Suppose that $P_1,\ldots,P_l$ are already hyperbolic while $P_{l+1}$ is not, where $l\ge 0$. By recursion, it is sufficient
to move $P_{l+1}$ off $U_0\cup V_0$ while leaving the points $P_1,\ldots,P_l$ hyperbolic. It is enough to consider the 
following two cases:

\noindent Case 1: $u_{l+1}=0$, $v_{l+1}\ne 0$, and

\noindent Case 2: $u_{l+1}=v_{l+1}=0$.

\noindent We claim that there exists an automorphism $\varphi\in\SAut(X)$ leaving $P_1,\ldots,P_l$ hyperbolic
such that  in Case 1 the point $\varphi\cdot P_{l+1}$ is hyperbolic as well, and  in Case~2 this point  satisfies the assumptions of Case 1.

In Case 1 we divide $P_1,\ldots,P_{l+1}$ into several disjoint pieces $M_0,\ldots,M_k$ according to different values of $v$ so
that $P_i\in M_j\Leftrightarrow v_i=c_j$, where $c_j\neq 0$. Assuming that~$M_0=\{P_{i_1}, \ldots, P_{i_r}, P_{l+1}\}$,
where $i_k\le l$ for all $k=1,\ldots,r$, we can choose an extra point $P'_{l+1}\in (V_{c_0}\cap X^1) \setminus U_0$.
Indeed, since $c_0=v_{l+1}\neq 0$, we have $V_{c_0}\cong Y^1$. We have $\dim Y^1\ge 2$,
hence $\dim (V_{c_0}\cap X^1) \setminus U_0=\dim Y^1\ge 2$.

By Lemma~\ref{stablevel} the subgroup $\Stab^v_{c_1,\ldots,c_{k}}\subseteq G_v$ acts
$(r+1)$-transitively on ${V_{c_0}\cap X^1}$. Therefore we can send the $(r+1)$-tuple $(P_{i_1},\ldots, P_{i_r},
P_{l+1})$ to \linebreak$(P_{i_1},\ldots, P_{i_r}, P'_{l+1})$ fixing the remaining points of $M_1\cup \ldots \cup M_{k}$. This confirms our claim in Case~1.

In Case 2 we have $P_{l+1}=(R_{l+1},0,0)\in X^1$. It follows from Lemma~\ref{regular} that
$R_{l+1}=\pi(P_{l+1})$ belongs to~$Y_{\reg}$ and $df(R_{l+1})\neq 0$ in the cotangent space $T^*_{R_{l+1}} Y$. The variety
$Y^1$ being flexible, there exists an LND $\partial_0\in\Der \Real[Y]$ such that $\partial_0(f)(R_{l+1})\neq 0$. Let
$q(v)=v(v-v_1)(v-v_2)\ldots(v-v_l)$ be a polynomial in~$\Real[v]$ and choose a set of generators $x_1,\ldots,x_s$ of the algebra $\Real[Y]$.
Then, as in (\ref{drob}), the derivation $\partial_0$ can be extended to $\partial_1\in \Der\, \Real[X]$ via
\begin{eqnarray*}
\partial_1 (x_i) &=& q(v)\partial_0 (x_i),\quad i=1,2,\ldots, s\;,\\
\partial_1 (v) &=&
0\;,\\
\partial_1 (u) &=& \frac{q(v)}v\, \partial_0 (f)\,.
\end{eqnarray*}
Due to our choice, $\partial_1 (u)(P_{l+1})\neq 0$. Hence the action of the associate one-parameter unipotent subgroup
$H(\partial_0,q)=\exp (t\partial_1)$ pushes the point $P_{l+1}$ out of~$U_0$. So the orbit $H(\partial_0,q)\cdot P_{l+1}$ meets the
hypersurface $U_0\subseteq X^1$ in finitely many points. Similarly, for every $j=1,2,\ldots,l$ the orbit~$H(\partial_0,q)\cdot P_j\not\subseteq U_0$ meets~$U_0$ in finitely many points. Let $\varphi=\exp(t_0
\partial_1)\in H(\partial_0,q)\subseteq G_v$. For a general value of $t_0\in\Real$ the image $\varphi\cdot P_j$ lies
outside $U_0$ for all ${j=1,2,\ldots, l+1}$. Since the group $H(\partial_0,q)$ preserves $v$, the points
$\varphi\cdot P_1,\ldots,\varphi\cdot P_l$ are still hyperbolic. Interchanging the roles of $u$ and $v$, we achieve that the assumptions of Case 1 are fulfilled for the new
collection~$\varphi\cdot P_1,\ldots,\varphi\cdot P_l,\varphi\cdot P_{l+1}$, as required.
\end{proof}

\section{Infinite transitivity on one connected component}

This section is devoted to the proof of Theorem~\ref{vtcc}. Recall that $Y$ is an affine variety defined over $\Real$, $f\in\Real[Y]$ is non-constant, and $X=\SP(Y,f)$. Let~$Y^1$ be a connected component of~$Y_{\reg}$, and let $X^1$ be a connected component
of~$\SP(Y^1,f)\cap X_{\reg}$.

\begin{lemma}\label{choicealpha}
Let $m$ be a positive integer and let $P_1,\ldots,P_m$ be $m$ points in $X^1$. 
There exist an automorphism $g\in\SAut(X)$ and a nonzero real number~$\alpha$ such that for each $i=1,\ldots,m$ the number $\alpha v(g\cdot P_i)$ is an interior point of
$f|_{Y^1}$.

Moreover, for any finite sets of real numbers $\mathcal U$ disjoint from $\{u(P_i)\}$ and $\mathcal V$ disjoint from~$\{v(P_i)\}$, the automorphism $g$ can be chosen in $\langle \Stab^u_{\mathcal U},\Stab^v_{\mathcal V} \rangle$.
\end{lemma}

\begin{proof}
Acting with $G_v$, we may assume that the $m$ points~$P_1,\ldots,P_m$ have pairwise distinct~$u$-coordinates. Acting further with $G_u$, we may assume that these points have also pairwise distinct values of their~$v$-coordinates. 

The proof depends on the behavior of $f$.
If $0$ is an interior point of~$f(Y^1)$, we let $g=\Id$ and choose $\alpha$ small enough. Then all $\alpha v(P_i)$ are close to~$0$, and thus are interior points of~$f(Y^1)$.

If $f|_{Y^1}$\label{case:1} is  unbounded  and $f|_{Y^1}>0$, we let $g=\Id$ and choose $\alpha$ great enough. All $\alpha v(P_i)$ are then large enough and are thus interior points of~$f(Y^1)$. In the case $f|_{Y_1}$ is  unbounded  and $f|_{Y_1}<0$, the same argument works.

It remains to consider a bounded function $f$, that is $\overline{f(Y^1)}=[a,b]$. This splits into two cases:
either $ab\neq 0$, or $ab=0$.

Case~1.
Without loss of generality we may suppose that $0<a<b$. Let $v_i=v(P_i)$ and~$u_i=u(P_i)$.
Consider the connected component $X^1$ of the suspension over~$Y^1$ such that all $v_i > 0$, and all $u_i > 0$. Let $v_{\max}$ and $v_{\min}$ be the maximal and minimal
values of~$v(P_1),\ldots,v(P_m)$.

If $v_{\max}/v_{min} \,<\, b/a$, we let $g=\Id$. Then for any $\alpha\in ] a/{v_{\min}}, \, b/{v_{\max}}[$, it is clear that all real numbers $\alpha v_i$
are interior points of $f(Y^1)$. 

Otherwise, if $v_{\max}/v_{min} \geqslant b/a$, we need a non-trivial automorphism $g$. 
Fix~$\varepsilon>0$. 
Note that, acting with the group $G_v$, any point $P\in\{P_1,P_2,\dots,P_m\}$ can be mapped to a point $P'$ such that $f(P')$ is very close to~$b$, while all the other points are fixed. Indeed, 
for a general~$\varepsilon_1<\varepsilon$, the real number $\frac{b-\varepsilon_1}{v(P)}\not\in \{u_1,\dots,u_m\}$. Let $R$ in $Y^1$ be such that $f(R)=b-\varepsilon_1$. We endow~$R$ with two extra coordinates $u=\frac{b-\varepsilon_1}{v(P)}, v=v(P)$ and get $P'=(R,\frac{b-\varepsilon_1}{v(P)}, v(P))\in X^1$. Let $\mathcal W=\{v(P_1),\dots,v(P_m)\}\setminus \{v(P)\}$, the point~$P$ can be mapped to $P'$ by an element $g_1$ of the  group~$\Stab^v_{\mathcal W}$. 
%
Thus for a general $\varepsilon_1<\varepsilon$, the automorphism $g_1$ of $X$ satisfies $g_1\cdot P_j = P_j$ for $P_j \ne P$, and~$f(g_1\cdot P) \,>\, b-\varepsilon$. 

Choose $P=P_i$ such that $v(P_i)=v_{\max}$, $P_i=(R_i,u_i,v_{\max})$. As described above, we map $P_i$ to $P_i'=(R',\frac{b-\varepsilon_1}{v_i}, v_i)$ such that $f(R')=b-\varepsilon_1$. Then, 
interchanging $u$ and $v$, and interchanging $a$ and $b$, there is an element of $G_u$ which maps $P_i'$ to $P_i''=(R'',\frac{b-\varepsilon_1}{v_i}, \frac{(a+\varepsilon_2)v_i}{b-\varepsilon_1})$ such that $f(R'')=a+\varepsilon_2$.
Note that $v(P_i'')<\frac{a+\varepsilon}{b-\varepsilon} v(P_i)$. 

If $v_{\max}/v_{\min}\geq b/a$ for the new set $P_1,\ldots,P_i'', \ldots, P_m$, we repeat this procedure. This process is finite since at each step the product $v(P_1)\ldots v(P_m)$ reduces by a factor at least~$\frac{a+\varepsilon}{b-\varepsilon}$. Finally, we get $m$ points $g\cdot P_1,\ldots,g\cdot P_m$ such that $v_{\max}/v_{min} \,<\, b/a$.

Case~2. One of $a$ and $b$ equals zero.
Using Lemma~\ref{avoid0}, we map the given $m$-tuple $P_1,\ldots,P_m$ in $X^1$ to points $P'_1,\ldots,P'_m \in X^1 \setminus (U_0\cup V_0)$. 
A sufficiently small $\alpha$ then fulfills the required condition.

To prove the second part of the lemma, we run the proof once again but we add an extra condition while performing the lift of an automorphism from $Y$ to $X$. Namely, for an automorphism in~$G_u$, we multiply the polynomial $q$ by $\prod_{u\in\mathcal U}q(u)$, and for an automorphism in~$G_v$, we multiply the polynomial $q$ by $\prod_{v\in\mathcal V}q(v)$.
\end{proof}

 \begin{proof}[Proof of Theorem~\ref{vtcc}]
Fix two $m$-tuples of distinct points $P_1,\ldots,P_m$ and $Q_1,\ldots,Q_m$ in $X^1$. By Lemma~\ref{choicealpha}, up to the action of $\SAut(X)$, there exists $\alpha\in\Real\setminus \{0\}$ such that $\alpha v(P_i)$ belongs to~$\Int f(Y^1)$ and $\alpha v(Q_i)\in \Int f(Y^1)$ 
for all $i=1,\ldots,m$.  
We denote by  $c_1,\ldots,c_k$ the distinct values of the $v$-coordinates of the given $2m$ points.   We split the set $\{P_1,\ldots,P_m,Q_1,\ldots,Q_m\}$ into~$k$ subsets according to the $v$-coordinate. For each $i$, the set
 $V_{c_i}\cap U_{\alpha}\cap X^1$ is infinite by Lemma~\ref{level}. In particular, for each~$i$, we get $\#\left(V_{c_i}\cap U_{\alpha}\cap X^1\right)\geq 2m$.
 By Lemma~\ref{stablevel},  there exists $g_i\in \Stab^v_{c_1,\ldots,\check c_i,\ldots,c_k}$ such that 
 $g_i\bigl(\{P_1,\ldots, P_m, Q_1, \ldots Q_m\}\cap V_{c_i}\bigr)\subset V_{c_i}\cap U_{\alpha}$ and $g_i$ fixes all the points belonging to $\cup_{j\neq i} V_{c_j}$.  Let us denote by $P'_1,\ldots,P'_m,Q'_1,\ldots,Q'_m\in U_{\alpha}$ the images by $g_k\circ\ldots \circ g_1$. Since the action of $\Stab^u_{\emptyset}=G_u$ on $U_{\alpha}$ is infinitely transitive by  Lemma~\ref{stablevel}, 
there exists a special automorphism mapping the  
 $m$-tuple $P'_1,\ldots,P'_m$ to $Q'_1,\ldots,Q'_m$. 
 \end{proof}
 
 \begin{lemma}
If $\SAut Y$ acts infinitely transitively on $Y^1$ where $\dim Y^1\geqslant 2$, then the flexibility of $Y^1$ implies the flexibility of any connected component $X^1$ of $\SP(Y^1,f)$, 
\end{lemma}

\begin{proof}

Clearly, $X^1$ is flexible if one point $P=(R,u,v)\in X^1$ is and if the group $\SAut(X)$ acts transitively on $X^1$.
Since the function $f\in\Real[Y^1]$ is non-constant,  $df(R)\neq 0$ at some point $R\in Y^1$ with $f(R)\neq 0$. Due to our
assumption~$Y^1$ is flexible. Hence there exist $n$ locally nilpotent derivations
$\partial_0^{(1)},\ldots,\partial_0^{(n)}\in\Der\Real[Y]$, where $n=\dim Y =\dim Y^1$, such that the corresponding vector fields
$\xi_1,\ldots, \xi_n$ span the tangent space $T_{R}Y$, i.e.
$$
\rk \left(%
\begin{array}{c}
  \xi_1(R) \\
  \vdots \\
  \xi_n(R) \\
\end{array}%
\right)=n\,.
$$
It follows that $\partial_0^{(i)} (f)(R)\ne 0$ for some index $i\in\{1,\ldots,n\}$.

Let now $P=(R,u_0,v_0)\in X^1$ be a point such that $\pi(P)=R$. Since $u_0v_0=f(R)\neq 0$, the point $P$ is
hyperbolic. Performing a lift as in~(\ref{drob}) with $q(v)=v$, we obtain $n$ LNDs
$$
\partial_1^{(1)},\ldots,\partial_1^{(n)}\in\Der\Real[X],\quad\mbox{ where}\quad
\partial_1^{(j)}=\partial_1^{(j)}(\partial_0^{(j)},v).
$$
Interchanging $u$ and $v$ and letting $j=i$, we get another LND 
$$
\partial_2^{(i)}=\partial_2^{(i)}(\partial_0^{(i)},u) \in\Der\Real[X].
$$
Let us show that the corresponding $n+1$ vector fields span the tangent space $T_{R}X$ at~$R$, as required. We can view
$\partial_1^{(1)},\ldots,\partial_1^{(n)},\partial_2^{(i)}$ as LNDs in $\Der\Real[Y][u,v]$ preserving the ideal~ $(uv-f)$, so
that the corresponding vector fields are tangent to the hypersurface 
$$
X=\{uv-f(R)=0\}\subseteq Y\times\Affspace^2\,.$$ The values  of these vector fields at the point $P'\in X$ yield an $(n+1)\times (n+2)$-matrix
$$
E=\left(%
\begin{array}{ccc}
  v_0\xi_1(R) & \partial_0^{(1)}(f)(R) & 0 \\
  \vdots & \vdots & \vdots \\
  v_0\xi_n(R) & \partial_0^{(n)}(f)(R) & 0 \\
   u_0\xi_i(R) & 0 & \partial_0^{(i)}(f)(R) \\
\end{array}%
\right).
$$
The first $n$ rows of $E$ are linearly independent, and the last one is independent from the preceding since
$\partial_0^{(i)}(f)(R)\neq 0$. Therefore $\rk (E)=n+1=\dim X$. So these locally nilpotent vector fields indeed span the
tangent space $T_{P}X$ at $P$, as claimed.
 \end{proof}
 
\section{Infinite transitivity on each connected component}
 
In this section we prove Theorem~\ref{main}. As above, $Y$ is an affine variety defined over $\Real$, $f\in\Real[Y]$ is non-constant, and $X=\SP(Y,f)$. Moreover, we assume that~$\SAut(Y)$ acts infinitely transitively on each connected component of $Y_{\reg}$.
Recall that~$v(P)$ denotes the $v$-coordinate of the point $P\in X \subseteq Y\times \Spec \Real[u,v]$.

\begin{lemma}\label{alldifferent}
Assume that $Y$ is a flexible variety. For every finite set of points $P_1, \ldots,P_m$ in~$X_{\reg}$, there exists an automorphism~$g\in\SAut(X)$ such that all $v(g\cdot P_i)$ are pairwise distinct, and all 
$u(g\cdot P_i)$ are pairwise distinct. 
\end{lemma}

\begin{proof}
We cannot use the starting argument of the proof of Lemma~\ref{choicealpha}, since several points of~$X$ may have the same projection in $Y$.

We denote the set of projections  $\pi(P_1), \ldots,\pi(P_m)$ by $\{ R_1,\ldots,R_{m'}\}$. Note that all~$R_i$ belong to~$Y_{\reg}$ by 
Lemma~\ref{regular}. 
Up to a special automorphism of $X$ we can assume  that all $f(R_i)$ are pairwise distinct. This is possible since $f$ is non-constant on each connected component, and 
 the action of $\SAut(Y)$ on $Y$ is infinite transitive on each connected component. Consider the images of~$P_1, \ldots,P_m$ under the projection $\rho\colon X \to \Spec \Real[u,v]$. If $\pi(P_i)\neq \pi(P_j)$, the projections $\rho(P_i)$ and~ $\rho(P_j)$ cannot coincide. Otherwise $f(R_i)=u_iv_i=u_jv_j=f(R_j)$. If $\pi(P_i)=\pi(P_j)$, we get also $\rho(P_i)\ne \rho(P_j)$ since the points $P_i$ and $P_j$ are distinct. 

Thanks to Lemma~\ref{avoid0}, keeping $\rho(P_i)\ne \rho(P_j)$ if $P_i\ne P_j$, we may assume that $u(P_i)\ne 0$ and $v(P_i)\ne 0$ for each $i$. 
We split the set $\{P_1, \ldots,P_m\}$ into several subsets $M_1\sqcup \ldots \sqcup M_k$ according to the $v$-coordinate. Let $c_i\in\Real$ be such that $M_i\subseteq V_{c_i}$. 
Using Lemma~\ref{stablevel}, we act by  an element of $\prod_{i=1}^k\Stab^v_{c_1,\ldots, \check c_i,\ldots,c_k}$ to get $m$ points with
pairwise distinct $u$-coordinates. Arguing likewise with $\Stab^u$-actions, we get $m$ points with
pairwise distinct $u$- and $v$-coordinates.
\end{proof}

\begin{proof}[Proof of Theorem~$\ref{main}$]
We denote by $s$ the number of connected components of~$Y_{\reg}$ and we suppose that the action of $\SAut(Y)$ on $Y_{\reg}$ is infinitely transitive on each connected component. Consider a suspension $X=\SP(Y,f)$ and decompose $X_{\reg}=X^1\sqcup\ldots\sqcup X^{s'}$  into connected components. Recall that over each connected component of $Y_{\reg}$ there is either one or two connected components of $X_{\reg}$. Fix some integers $m'_1,m'_2,\ldots,m'_{s'}$ such that $\sum m'_j = m$ and choose two $(m'_1+\ldots+m'_{s'})$-tuples $\mathcal{P}=\{P_1,\ldots,P_m\}$ and $\mathcal{Q}=\{Q_1,\ldots,Q_m\}$ in $X$ such that for each $j$,  the component $X^j$ contains exactly $m'_j$ points of $\mathcal{P}$ and~$m'_j$ points of $\mathcal{Q}$. Let~$\mathcal{S}=\mathcal{P}\cup\mathcal{Q}=\{S_1,\dots S_{2m}\}$.

By Lemma~\ref{alldifferent} applied to $\mathcal{S}$, we may suppose that the values of the $v$-coordinates are pairwise distinct and that the values of the $u$-coordinates are also pairwise distinct. 

We want to choose an $s'$-tuple of values $\alpha=(\alpha_1,\ldots,\alpha_{s'})$ such that for all $S_i\in X^j$ the number $\alpha_jv(S_i)$ is an interior point of $f|_{X^j}$.
To this end, we repeatedly apply Lemma~\ref{choicealpha} proceeding on one connected component at each step.
Notice that we need to preserve the condition that the values of the $v$-coordinates are pairwise distinct and that the values of the $u$-coordinates are also pairwise distinct. 
At $j$th step, we let $\mathcal{U}=\{u(S_i),S_i\not\in X^j\}$ and~$\mathcal{V}=\{v(S_i),S_i\not\in X^j\}$ and we use $g\in \langle \Stab^u_{\mathcal U},\Stab^v_{\mathcal V} \rangle$ given by Lemma~\ref{choicealpha}.

Such a choice of $g$ provides that $g\cdot S_i = S_i$ for $S_i\not \in X^j$. To control the condition that all~$u$-values are pairwise distinct and all $v$-values are pairwise distinct for 
$S_i\in X^j$, we require the following. For each one-parameter subgroup $h(t)$ acting in the course of the proof of Lemma~\ref{choicealpha} (recall that it is non-trivial only for $S_i\in X^j$),
the conditions on $t$ are
\begin{eqnarray*}
&v(h(t)\cdot S_i)\not \in \mathcal V, \, u(h(t)\cdot S_i)\not \in \mathcal U \; \mbox{for } S_i \in X^j\, ;
\\
&v(h(t)\cdot S_{i_1})\neq v(h(t)\cdot S_{i_2}), \; u(h(t)\cdot S_{i_1})\neq u (h(t)\cdot S_{i_2}) \; \\
&\mbox{for distinct }  S_{i_1} , S_{i_2}\in X^j\;.
\end{eqnarray*}
This is true for generic $t$. At each step we get an $\alpha_j$ which fits for all $S_i\in X^j$. We may choose the $\alpha_i$ pairwise distinct. At the end, we get a collection $\alpha=(\alpha_1,\dots,\alpha_{s'})$, as required.

We construct an automorphism $g_0$ mapping $\mathcal S$ to $(X^1\cap U_{\alpha_1})\sqcup (X^2\cap U_{\alpha_2}) \sqcup \ldots \sqcup (X^{s'}\cap U_{\alpha_{s'}})$ as the product of $2m$ automorphisms, each of them fixing all the points but one.
Since all~$v$-values are pairwise distinct, to map $S_i$, we let $q(v)=\beta v\prod _{k\ne i} (v-v(S_k))$
where $\beta$ satisfies $q(v(S_i))=1$.
If $S_i\in X^j$, using the lift defined by $q$ (see~(\ref{drob})), we map $S_i$ to $X^j\cap U_{\alpha_j}$. 
Notice that for $\{g_0\cdot S_1,\ldots,g_0\cdot S_{2m}\}$, the $u$-values are no longer pairwise distinct. 

To map $g_0\cdot P_1,\ldots,g_0\cdot P_m$ onto $g_0\cdot Q_1,\ldots,g_0\cdot Q_m$, we  use, for each $i$, the infinite transitivity of the group $\Stab^u_{\alpha_1,\ldots,\check \alpha_i,\ldots,\alpha_{s'}}$, multiplying the corresponding LNDs on $\Real[Y]$ by the polynomial~ $q(u)=\gamma u(u-\alpha_1)\ldots\check{\overbrace{(u-\alpha_i)}}\ldots(u-\alpha_{s'})$ where $\gamma$ is such that $q(\alpha_i)=1$. 
In this way, for each~$i$, we fix the points lying off the $i$th connected component. Finally we get an automorphism~$g$ which maps $g_0\cdot P_1,\ldots,g_0\cdot P_m$ to $g_0\cdot Q_1,\ldots,g_0\cdot Q_m$. Hence, $g_0^{-1}gg_0$ maps the $m$-tuple~$P_1,\ldots,P_m$ to $Q_1,\ldots,Q_m$, as required. 
\end{proof}

\section*{Acknowledgements}
The first author benefited from the support of the "EADS Foundation Chair in Mathematics", RFBR grant 09-01-00648a, and Russian-French Poncelet Laboratory (UMI 2615 of~CNRS).

\end{document}